\documentclass[10pt]{article}
\usepackage{amsfonts}
\usepackage{amssymb}
\usepackage{graphics}

\newtheorem{theorem}{\sc Theorem}[section]

\newtheorem{prop}{\sc Proposition}[section]
\newtheorem{remark}{\sc Remark  }

\newcommand{\bce}{\begin{center}}
\newcommand{\ece}{\end{center}}

\begin{document}

\title{Detecting the limit cycles for a class of Hamiltonian systems under
thirteen-order perturbed terms }

\author{Gheorghe Tigan \thanks{%
Department of Mathematics, "Politehnica" University of Timisoara,
Timisoara, Timis, Romania; email: gtigan73@yahoo.com} }
\date{}
\maketitle

\begin{abstract}

It this paper we study a class of perturbed Hamiltonian systems
under perturbations of thirteen order in order to detect the
number of limit cycles which bifurcate from some periodic orbits
of the unperturbed Hamiltonian system. The system has been
previously studied in \cite{tig2}, \cite{tig3}, \cite{tig4}. We
observe in the present work that the system under perturbations of
thirteen order can have more limit cycles than under perturbations
of five \cite{tig3}, respectively nine \cite{tig4} order but we
have not identified more limit cycles than under perturbations of
seven \cite{tig2} order.

\end{abstract}
\emph{Key words:} Hamiltonian systems, limit cycles, Abelian
integral.\\
AMS 2000: 34C07, 37G15 \\

\section{Introduction}

Consider the following perturbed Hamiltonian system:
\begin{equation}\label{sis3n12}
\left\{
\begin{array}{c}

\dot{x}=4y\left( abx^{2}-by^{2}+1\right) +\varepsilon
x(ux^{n}+vy^{n}-b\frac{\beta +1}{\mu +1}x^{\mu }y^{\beta
}-ux^{2}-\lambda ), \\

\dot{y}= 4x\left( ax^{2}-aby^{2}-1\right) +\varepsilon
y(ux^{n}+vy^{n}+bx^{\mu }y^{\beta }-vy^{2}-\lambda ) \ \ \ \ \ \
\end{array}
\right.
\end{equation}

where $\mu +\beta =n,0<a<b<1,0<\varepsilon \ll 1,$ $u,v,\lambda $
are the real parameters and $n=2k, k$ an integer positive.
Studying the existence, number and distribution of limit cycles in
a system of polynomial differential equations is an open problem
even for a 2-degree polynomial and it is known as the Hilbert`s
16th problem.
\par The system of differential equations:

\begin{equation}\label{li}
\left\{
\begin{array}{c}
\dot{x}=\allowbreak y\left( 1+x^{2}-ay^{2}\right) +\varepsilon
x(mx^{2}+ny^{2}-\lambda ) \ \ \\
\dot{y}=-\allowbreak x\left( 1-cx^{2}+y^{2}\right) +\varepsilon
y(mx^{2}+ny^{2}-\lambda )%
\end{array}%
\right.
\end{equation}

\noindent has been discussed in \cite{Li}, \cite{tig1} and the
system

\begin{equation}
\left\{
\begin{array}{c}
\dot{x}=\allowbreak y\left( 1-cy^{2}\right) +\varepsilon
x(mx^{2}+ny^{2}-\lambda ) \ \ \\
\dot{y}=-\allowbreak x\left( 1-ax^{2}\right) +\varepsilon
y(mx^{2}+ny^{2}-\lambda )%
\end{array}%
\right.
\end{equation}
 \noindent in \cite{Li4}. It has been  shown  that each of the
two systems can have at least 11 limit cycles.

The higher order perturbations of the system (\ref{li}),  have
been recently studied in \cite{Cao}, \cite{Tang}.

The following result is reported in \cite{Cao}.

\begin{theorem}\label{teorema1} Consider the perturbed Hamiltonian system

\begin{equation}\label{sistem1}
\dot{x}=-\frac{\partial H}{\partial y}+P(x,y,\alpha ),\dot{y}=%
\frac{\partial H}{\partial x}+Q(x,y,\alpha )
\end{equation}

Assume that $P(x,y,0)=Q(x,y,0)=0$, and the unperturbed system
exhibits a center. Denote by  $\Gamma ^{h}$ the closed curves
$H(x,y)=h$ surrounding this center, and by $\Gamma ^{h}(D)$ the
subset bounded by $\Gamma ^{h}$ and containing the center. Suppose
that as $h$ increases the diameter of the set $\Gamma ^{h}(D)$
increases.
 If there exists $h_{0}$ such that function

\begin{equation}
A(h)=\int\limits_{\Gamma ^{h}(D)}[P_{x\alpha }^{\prime \prime
}(x,y,0)+Q_{y\alpha }^{\prime \prime }(x,y,0)]dxdy
\end{equation}

\noindent satisfies $A(h_{0})=0,A^{\prime }(h_{0})\neq 0,\alpha
A^{\prime }(h_{0})<0, respectively >0,$ then the perturbed system
(\ref{sistem1}) has only one stable, respectively unstable limit
cycle nearby $\Gamma ^{h_{0}}$, for $\alpha $ very small. If
$\Gamma ^{h}$ shrinks  as $h$ increases, the stability of the
limit cycle is opposite to the above cases. If $A(h)\neq 0,$ then
the system (\ref{sistem1}) has no limit cycle.
\end{theorem}
The integral $A(h)$ is called the \emph{Abelian integral}
\cite{Blows}. If the Hamiltonian system (\ref{sistem1}) is
perturbed in the form:

\begin{equation}\label{sistem2}
\left\{
\begin{array}{c}
\dot{x}(t)=-\frac{\partial H}{\partial y}+\varepsilon
x(p(x,y)-\lambda ),
\\
\dot{y}(t)=\frac{\partial H}{\partial x}+\varepsilon
y(q(x,y)-\lambda ), \ \ \
\end{array}%
\right.
\end{equation}

\noindent where $p(0,0)=q(0,0)=0,$ then, by the above Theorem
\ref{teorema1}, from $A(h)=0$, we get:

\begin{equation}\label{detec1}
\lambda =\lambda (h)=\frac{\int\limits_{\Gamma^{h}(D)}f(x,y)dxdy}{%
2\int\limits_{\Gamma^{h}(D)}dxdy},
\end{equation}

\noindent with $f(x,y)=xp_{x}^{\prime }+yp_{y}^{\prime }+p+q.$
This function $\lambda (h)$ is called \emph{the detection
function} of the system (\ref{sistem2}).

Using the detection function $\lambda (h)$ we get by Theorem
\ref{teorema1}  the following result :

\begin{prop} \label{prop1}  a) If $(h_{0},\lambda (h_{0}))$ is an
intersecting point of a line $\lambda =\lambda _{0}$ and the detection curve $%
\lambda =\lambda (h)$, with  $\lambda ^{\prime }(h_{0})>0 (<0)$,
then the system (\ref{sistem2}) has only one stable (unstable)
limit cycle nearby $\Gamma ^{h_{0}}$;  b) If the line $\lambda
=\lambda _{0}$ and the detection curve $\lambda =\lambda (h)$ do
not intersect each other, then the
system (\ref{sistem2}) has no limit cycle.  If the $%
\Gamma ^{h}(D)$ shrinks  as h increases, the stability of the
limit cycle is opposite to the above cases.
\end{prop}

This paper is organized as follows. In Section 2, we recall the
results regarding the unperturbed system. In Section 3, using the
analytical expressions of the detection functions and a Computer
Algebra System we numerically compute the detection functions.
Based on these data we can illustrate the distribution of the
limit cycles.

\section{The behavior of the unperturbed system}

The unperturbed system (\ref{sistem4}) of the system
(\ref{sis3n12}) is:

\bigskip
\begin{equation}\label{sistem4}
\left\{
\begin{array}{c}
\dot{x}=\allowbreak 4y\left( -by^{2}+abx^{2}+1\right) \\
\dot{y}=\allowbreak 4x\left( ax^{2}-aby^{2}-1\right) \ \
\end{array}%
\right.
\end{equation}
The Hamilton function defining this system is:

\begin{equation}\label{hamil2}
H\left( x,y\right) =-(ax^{4}+by^{4})+2abx^{2}y^{2}+2\left(
x^{2}+y^{2}\right) =h
\end{equation}

\noindent The function $H$  has nine finite singular points:

$O(0,0),A_{1}\left( \sqrt{\frac{1+a}{a\left( 1-ba\right) }},\frac{1}{b-b^{2}a%
}\sqrt{b\left( 1-ba\right) \left( 1+b\right) }\right) ,A_{2}\left( \sqrt{%
\frac{1+a}{a\left( 1-ba\right) }},\frac{-1}{b-b^{2}a}\sqrt{b\left(
1-ba\right) \left( 1+b\right) }\right) ,$

$A_{3}\left( -\sqrt{\frac{1+a}{a\left( 1-ba\right) }},\frac{1}{b-b^{2}a}%
\sqrt{b\left( 1-ba\right) \left( 1+b\right) }\right) ,A_{4}\left( -\sqrt{%
\frac{1+a}{a\left( 1-ba\right) }},\frac{-1}{b-b^{2}a}\sqrt{b\left(
1-ba\right) \left( 1+b\right) }\right) ,$  \\
$B_{1}(0,\sqrt{\frac{1}{b}}),B_{2}(0,-\sqrt{\frac{1}{b}}),C_{1}(\sqrt{\frac{1}{a}},0)$ and $C_{2}(-%
\sqrt{\frac{1}{a}},0)$ .

Computing the eigenvalues of the associated linear system at each singular point we get that the points $%
O, A_{1},A_{2},A_{3},A_{4}$ are centers, while
$B_{1},B_{2},C_{1},C_{2}$ are  saddle points.

The values of the Hamilton function $H$ at the singular points are respectively: $H\left( A_{i}\right) =\allowbreak \frac{2ba+b+a}{ba\left( 1-ba\right) }%
,i=1-4,$ $H\left( B_{k}\right) =\allowbreak \frac{1}{b},H\left(
C_{k}\right) =\allowbreak \frac{1}{a},k=1,2.$ Because $0<a<b<1$ we
get that: $H\left( O\right) <H\left( B_{1}\right) <H\left(
C_{1}\right) <H\left( A_{1}\right). $

In polar coordinates, $x=r\cos \theta ,y=r\sin \theta ,$ the
system (\ref{sistem4}) becomes:

\begin{equation}\label{sistem5}
\dot{r}=-r^{3}p^{\prime }(\theta
),\,\,\,\dot{\theta}=-1+r^{2}p(\theta )
\end{equation}

\noindent and the Hamilton function (\ref{hamil2}) gets:
\begin{equation}\label{hamil3}
H(r,\theta )=-r^{4}p(\theta )+2r^{2}=h,
\end{equation}

\noindent where
\begin{equation}\label{pe}
p(\theta )=a\cos ^{4}\theta +b\sin ^{4}\theta -2ab\cos ^{2}\theta
\sin ^{2}\theta .
\end{equation}

\begin{remark} The equilibrium points $A_{1},A_{2},A_{3},A_{4}$ lie on the lines $%
d_{\pm }:$ $\theta =\pm \arccos \sqrt{\frac{b+ba}{a+b+2ab}}$
\end{remark}

\begin{figure}[h]\bce
\includegraphics{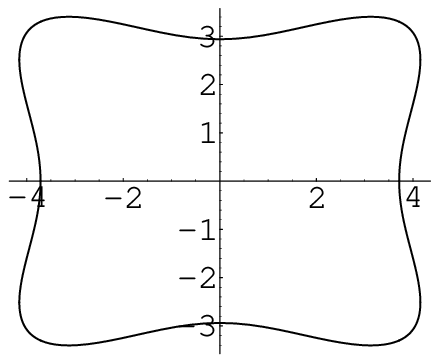}
\includegraphics{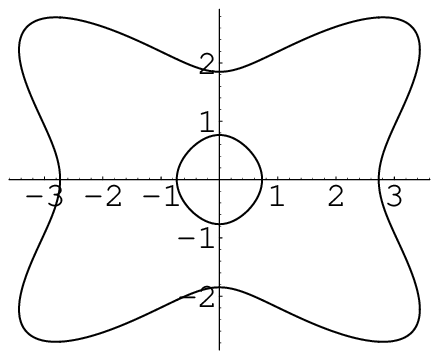}
\caption{Orbit of type a) $L_{1}$ (left)  b) $L_{2}$ and $L_{1}$
(right) } \label{neper31} \ece
\end{figure}

\begin{figure}[h]\bce
\includegraphics{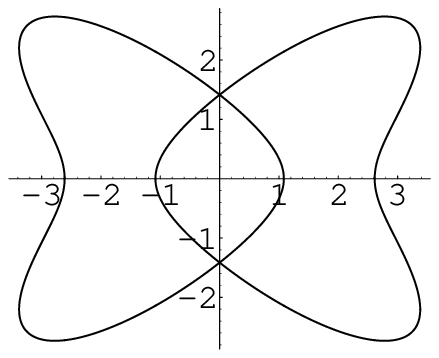}
\includegraphics{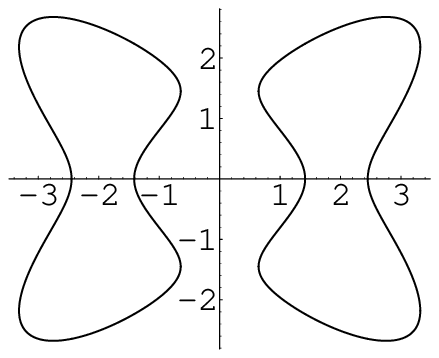}
\caption{ a) Four heteroclinic orbits connecting two critical
points $B_{1}$, $B_{2}$ (left) b) Two orbits of type $L_{3}$
(right) } \label{neper33} \ece
\end{figure}

\begin{figure}[h]\bce
\includegraphics{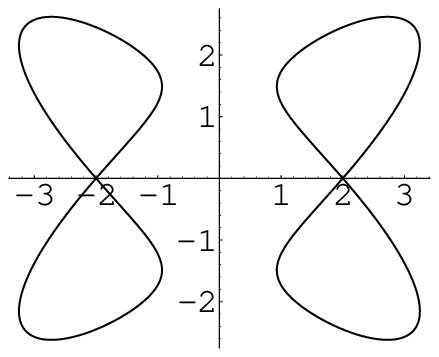}
\includegraphics{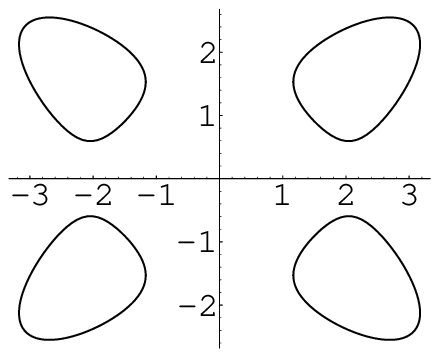}
\caption{a) Four homoclinic orbits connecting the critical points
$C_{1}$ and $C_{2}$ (left) b) Four orbits of type $L_{4}$ (right)
} \label{neper35} \ece
\end{figure}

\begin{theorem}\label{teorema2}

As $h$ varies on the real line, the level curves $H(x,y)=h$
 can be classified, as follows \cite{tig2}:
\begin{enumerate}
\item $\Gamma _{1}^{h}:-\infty <h<0,$ which corresponds to an
orbit that surrounds all critical points (Fig.\ref{neper31}a).

\item $\Gamma _{2}^{h}\cup \Gamma _{1}^{h}:0\leq h<\frac{1}{b},$
corresponding to an orbit $\left( \Gamma _{2}^{h}\right) $ that
surrounds only the point $O$ and a curve of type $\left( \Gamma
_{1}^{h}\right) ,$ (Fig.\ref{neper31}b,a).

\item $\Gamma _{3}^{h}:\frac{1}{b}<h<\frac{1}{a},$ which
corresponds to two symmetric orbits that do not cross the Oy axis,
but encircle the rest of the critical points. If $h= \frac{1}{b}$
we get four heteroclinic orbits connecting the critical points
$B_{1}$ and $B_{2}$ (Fig.\ref{neper33}b,a).

\item $\Gamma _{4}^{h}:$ \ $\frac{1}{a}<h<\frac{2ba+b+a}{ba\left( 1-ba\right) }%
, $ corresponding  to four orbits that surround respectively the points $%
A_{i},i=1-4$. If $h= \frac{1}{a}$ we have four homoclinic orbits,
namely, two homoclinic to  $C_{1}$ and two homoclinic to $C_{2}$
(Fig.\ref{neper35}b,a). Note that as $h$ increases, the curves
$\Gamma _{1}^{h},\Gamma _{3}^{h}$ and $\Gamma _{4}^{h}$ shrink,
while $\Gamma _{2}^{h}$ extends.
\end{enumerate}
\end{theorem}

\bigskip

\section{Numerical explorations}

\begin{figure}[h]\bce
\includegraphics{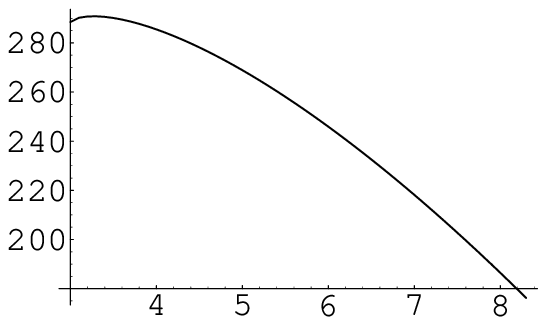}
\includegraphics{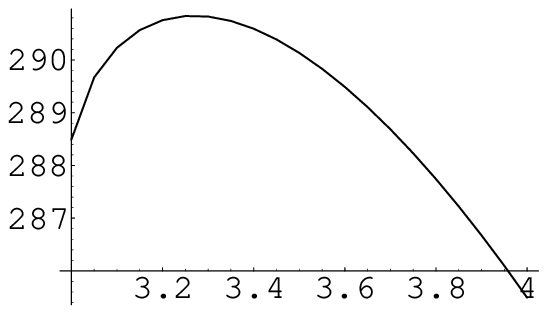}
\caption{Detection curve $\lambda_{4}$  of the system , for \ \
$a=\frac{1}{3}, b=\frac{1}{2}, n=12, u =0.007$ and $v=-0.028$. }
\label{detfig1} \ece
\end{figure}

\begin{figure}[h]\bce
\includegraphics{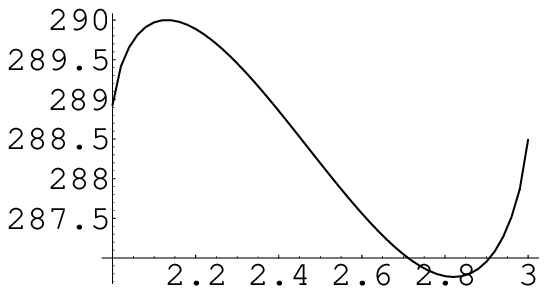}
\includegraphics{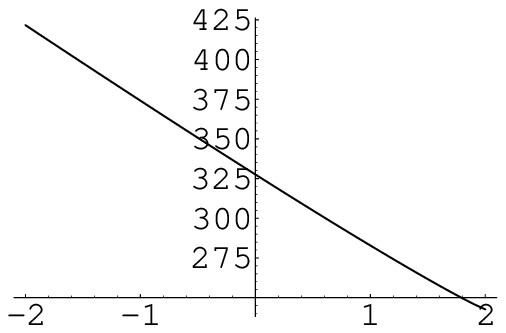}
\caption{Detection curves $\lambda_{3}$ (left), $\lambda_{1}$
(right) of the system (\ref{sis3n12}), for \ $a=\frac{1}{3},
b=\frac{1}{2}, n=12, u =0.007$ and $v=-0.028$ } \label{detfig2}
\ece
\end{figure}

In this section we numerically compute the detection curves and
point out the distribution of the limit cycles. The four detection
functions can be computed numerically, and for a given $h$, they
depend on $u$ and $v,$ (see tables 1-4). On the other hand, for two given values of $u$ and $%
v,$ the detection curves can be plotted on the $\left( h,\lambda \right) $%
-plane, as illustrated in Fig.\ref{detfig1}, \ref{detfig2}. By the
Proposition \ref{prop1} and the detection function graphs, the
existence, number and distribution of limit cycles can then be
illustrated. We consider here the case $n=12$.

From (\ref{hamil3}), we get
\begin{equation}\label{erurile}
r_{1,2}=r_{\pm }^{2}(\theta ,h)=\frac{1\pm \sqrt{1-hp(\theta
)}}{p(\theta )}
\end{equation}

and from $\dot{\theta}=-1+r^{2}p(\theta )=0$, we have: \\
\\
$\theta _{1}(h)=\frac{1}{2}\arccos \left[ \left( b-a+2\sqrt{%
a^{2}b^{2}-ab+\left( a+b+2ab\right) h^{-1}}\right) /\left( a+b+2ab\right) %
\right] ,$ \\
\\
$\theta _{2}(h)=\frac{1}{2}\arccos \left[ \left( b-a-2\sqrt{%
a^{2}b^{2}-ab+\left( a+b+2ab\right) h^{-1}}\right) /\left( a+b+2ab\right) %
\right] .$ \\
\\
Using the perturbation terms \\
\\
$p\left( x,y\right) =x\left( ux^{n}+vy^{n}-b\frac{\beta +1}{\mu
+1}x^{\mu }y^{\beta }-ux^{2}\right) ,q\left( x,y\right) =y\left(
ux^{n}+vy^{n}+bx^{\mu }y^{\beta }-vy^{2}\right) $ \\

we have $\frac{\partial ^{2}p\left( x,y\right) }{\partial
x\partial \varepsilon }+\frac{\partial ^{2}q\left( x,y\right)
}{\partial y\partial \varepsilon }=\left( 2+n\right) \left(
ux^{n}+vy^{n}\right) -3(ux^{2}+vy^{2})-2\lambda .$ \\
\\
Hence, the four detection functions, corresponding to the four
closed curves $\Gamma _{j}^{h},$ $j=1-4,$ for the above
perturbations are as follows: \\
\\
\begin{equation}\label{detec2}
\lambda _{j}(h)=\frac{\int\limits_{\Gamma _{j}^{h}(D)}\left[
(n+2)\left( ux^{n}+vy^{n}\right) -3(ux^{2}+vy^{2})\right]
dxdy}{2\int\limits_{\Gamma _{j}^{h}(D)}dxdy},j=1-4
\end{equation}

For $a=1/3, b=1/2$, and $n=12,$ we get by (\ref{detec2}) the next
four detection functions, in polar coordinates:

\bigskip
\begin{equation}\label{lam1}
\lambda _{1}(h)=\frac{\int\limits_{0}^{2\pi }\left(
r_{1}^{7}\left( \theta ,h\right) g\left( \theta \right)
-\frac{3}{4}r_{1}^{2}\left( \theta ,h\right) g_{1}\left( \theta
\right) \right) d\theta }{\int\limits_{0}^{2\pi }r_{1}\left(
\theta ,h\right) d\theta },\,\,-\infty <h<2,
\end{equation}

\begin{equation}\label{lam2}
\lambda _{2}(h)=\frac{\int\limits_{0}^{2\pi }\left(
r_{2}^{7}\left( \theta ,h\right) g\left( \theta \right)
-\frac{3}{4}r_{2}^{2}\left( \theta ,h\right) g_{1}\left( \theta
\right) \right) d\theta }{\int\limits_{0}^{2\pi }r_{2}\left(
\theta ,h\right) d\theta },\,\,0<h<2,
\end{equation}

\begin{equation}\label{lam3}
\lambda _{3}(h)=\frac{\int\limits_{-\theta _{2}\left( h\right)
}^{\theta _{2}\left( h\right) }\left[ \left( r_{1}^{7}\left(
\theta ,h\right)
-r_{2}^{7}\left( \theta ,h\right) \right) g\left( \theta \right) -\frac{3}{4}%
\left( r_{1}^{2}\left( \theta ,h\right) -r_{2}^{2}\left( \theta
,h\right) \right) g_{1}\left( \theta \right) \right] d\theta
}{\int\limits_{-\theta _{2}\left( h\right) }^{\theta _{2}\left(
h\right) }\left( r_{1}\left( \theta ,h\right) -r_{2}\left( \theta
,h\right) \right) d\theta },\,\, 2<h<3,
\end{equation}

\begin{equation}\label{lam4}
\lambda _{4}(h)=\frac{\int\limits_{\theta _{1}\left( h\right)
}^{\theta _{2}\left( h\right) }\left[ \left( r_{1}^{7}\left(
\theta ,h\right)
-r_{2}^{7}\left( \theta ,h\right) \right) g\left( \theta \right) -\frac{3}{4}%
\left( r_{1}^{2}\left( \theta ,h\right) -r_{2}^{2}\left( \theta
,h\right) \right) g_{1}\left( \theta \right) \right] d\theta
}{\int\limits_{\theta _{1}\left( h\right) }^{\theta _{2}\left(
h\right) }\left( r_{1}\left( \theta ,h\right) -r_{2}\left( \theta
,h\right) \right) d\theta },\,\, 3<h< 8.2,
\end{equation}

\bigskip

\noindent where $g\left( \theta \right) =u\cos ^{12}\theta +v\sin
^{12}\theta
,g_{1}\left( \theta \right) =u\cos ^{2}\theta +v\sin ^{2}\theta $ and $%
r_{1,2}\left( \theta ,h\right) =r_{\pm }^{2}\left( \theta
,h\right).$

\bigskip

Using the expressions (\ref{lam1})-(\ref{lam4}) and a Computer
Algebra System we find the values of the detection functions
$\lambda_{i}(h),i=1-4$, recorded in the following tables (1-4).
Denote by $\rho=10^{4}u$ and $\omega=10^{4}v$
\bigskip
\begin{center}
Table 1 \\
Values of the detection function $\lambda_{1}(h)$, for
$a=1/3, b=1/2, n=12.$ \\
\begin{tabular}{|c|c|c|c|c|c|}
  \hline
    $h$ & $\lambda_{1}(h)$ & $h$ & $\lambda_{1}(h)$ & $h$ & $\lambda_{1}(h)$ \\ \hline
      -2 & 4.933$\rho$+ 1.373$\omega$ &  -1.9 &
      4.862$\rho$+ 1.352$\omega$ &  -1.8 &
      4.792$\rho$+ 1.332$\omega$ \\ \hline

       -1.7 & 4.722$\rho$+ 1.311$\omega$ &  -1.6 &
      4.653$\rho$+ 1.291$\omega$ &  -1.5 &
      4.584$\rho$+ 1.271$\omega$ \\ \hline

      -1.4 & 4.516$\rho$+ 1.251$\omega$ &  -1.3 &
      4.448$\rho$+ 1.231$\omega$ &  -1.2 &
      4.381$\rho$+ 1.212$\omega$ \\ \hline

       -1.1 & 4.315$\rho$+ 1.192$\omega$ &  -1. &
      4.249$\rho$+ 1.173$\omega$ &  -0.9&
      4.183$\rho$+ 1.154$\omega$ \\ \hline

       -0.8 & 4.118$\rho$+ 1.135$\omega$ &  -0.7 &
      4.054$\rho$+ 1.116$\omega$ &  -0.6&
      3.990$\rho$+ 1.098$\omega$ \\ \hline

       -0.5 & 3.926$\rho$+ 1.079$\omega$ &  -0.4 &
      3.863$\rho$+ 1.061$\omega$ &  -0.3 &
      3.801$\rho$+ 1.043$\omega$ \\ \hline

      -0.2 & 3.739$\rho$+ 1.025$\omega$ &  -0.1 &
      3.678$\rho$+ 1.008$\omega$ &  0 &
      3.618$\rho$+ 0.9906$\omega$ \\ \hline

        0.1 & 3.558$\rho$+ 0.9733$\omega$ &  0.2 &
      3.498$\rho$+ 0.9562$\omega$ &  0.3 &
      3.439$\rho$+ 0.9392$\omega$ \\ \hline

       0.4 & 3.381$\rho$+ 0.9224$\omega$ &  0.5 &
      3.323$\rho$+ 0.9058$\omega$ &  0.6 &
      3.266$\rho$+ 0.8895$\omega$ \\ \hline

       0.7 & 3.210$\rho$+ 0.8733$\omega$&  0.8 &
      3.154$\rho$+ 0.8573$\omega$ &  0.9 &
      3.099$\rho$+ 0.8415$\omega$ \\ \hline

       1. & 3.044$\rho$+ 0.8259$\omega$ &  1.1 &
      2.990$\rho$+ 0.8105$\omega$ &  1.2 &
      2.937$\rho$+ 0.7954$\omega$ \\ \hline

      1.3 & 2.885$\rho$+ 0.7805$\omega$ &  1.4 &
      2.834$\rho$+ 0.7658$\omega$ &  1.5 &
      2.783$\rho$+ 0.7514$\omega$ \\ \hline

       1.6 & 2.734$\rho$+ 0.7373$\omega$&  1.7 &
      2.685$\rho$+ 0.7235$\omega$ &  1.8 &
      2.638$\rho$+ 0.7102$\omega$ \\ \hline

       1.9 & 2.593$\rho$+ 0.6974$\omega$ &  2. &
      2.553$\rho$+ 0.6859$\omega$ \\ \hline

\end{tabular}
\end{center}

\begin{center}
Table 2 \\
Values of the detection function $\lambda_{2}(h)$, for
$a=1/3, b=1/2, n=12.$ \\
\begin{tabular}{|c|c|c|c|}
  \hline
    $h$ & $\lambda_2(h)$ & $h$ & $\lambda_2(h)$  \\ \hline
0.01 &  -0.001875 u-0.001876 v &  0.11 &  -0.020732 u-0.02083 v \\
\hline
0.21 &  -0.039766 u-0.0401393 v & 0.31 & -0.058969 u-0.059813 v \\
\hline
0.41 &  -0.078305 u-0.079837 v & 0.51 & -0.097701 u-0.10015 v \\
\hline
0.61 &  -0.117021 u-0.120612 v & 0.71 & -0.136046 u-0.140967 v \\
\hline
0.81 &  -0.154444 u-0.160774 v & 0.91 & -0.171734 u-0.179321 v \\
\hline
1.01 &  -0.187236 u-0.195488 v & 1.11 & -0.200019 u-0.207541 v \\
\hline
1.21 & -0.208827 u-0.212808 v & 1.31 & -0.211988 u-0.207141 v \\
\hline
1.41 & -0.207302 u-0.184002 v & 1.51 & -0.19189 u-0.132774 v \\
\hline
1.61 &  -0.162016 u-0.0354102 v & 1.71 & -0.11285 u+0.141107 v \\
\hline
1.81 & -0.038213 u+0.465149 v & 1.91 & 0.069632 u+1.1178 v \\
\hline

\end{tabular}
\end{center}
\bigskip
\begin{center}
Table 3 \\
Values of the detection function $\lambda_{3}(h)$, for
$a=1/3, b=1/2, n=12.$ \\
\begin{tabular}{|c|c|c|c|c|c|}
  \hline
    $h$ & $\lambda_3(h)$ & $h$ & $\lambda_3(h)$ & $h$ & $\lambda_3(h)$ \\ \hline

    2 & 3.0408$\rho$ + 0.8167$\omega$ &  2.02 &  3.0446$\rho$ +
    0.8175$\omega$ & 2.04& 3.0459$\rho$ + 0.8177$\omega$ \\ \hline

    2.06 &   3.0463$\rho$ + 0.8177$\omega$ &  2.08 &
    3.0462$\rho$ + 0.8175$\omega$ &  2.1 &   3.0457$\rho$ + 0.8172$\omega$ \\ \hline

    2.12 &
    3.0449$\rho$ +0.8168$\omega$ &  2.14 &   3.0438$\rho$ + 0.8164$\omega$ &  2.16 &
    3.0425$\rho$ + 0.8159$\omega$ \\ \hline

    2.18 &   3.0411$\rho$ + 0.8154$\omega$ &  2.2 &
    3.0396$\rho$ + 0.8148$\omega$ &  2.22 &   3.0379$\rho$ + 0.8142$\omega$ \\ \hline

    2.24 &
    3.0361$\rho$ + 0.8136$\omega$ &  2.26 &   3.0342$\rho$ + 0.8130$\omega$ &  2.28 &
    3.0323$\rho$ + 0.8123$\omega$ \\ \hline

     2.3 &   3.0303$\rho$ + 0.8116$\omega$ &  2.32 &
    3.0283$\rho$ + 0.8110$\omega$ &  2.34 &   3.0262$\rho$ + 0.8103$\omega$ \\ \hline

     2.36 &
    3.0240$\rho$ + 0.8096$\omega$ &  2.38 &   3.0219$\rho$ + 0.8089$\omega$ &  2.4 &
    3.0197$\rho$ + 0.8082$\omega$ \\ \hline

    2.42 &   3.0176$\rho$ + 0.8075$\omega$ &
    2.44 & 3.0154$\rho$ + 0.8068$\omega$ &  2.46 &   3.0132$\rho$ +
0.8061$\omega$ \\ \hline
 2.48 &
    3.0110$\rho$ + 0.8054$\omega$ &  2.5 &   3.0089$\rho$ + 0.8047$\omega$ &  2.52 &
    3.0068$\rho$ + 0.8040$\omega$ \\ \hline

     2.54 &   3.0047$\rho$ + 0.8033$\omega$ &  2.56 &
    3.0026$\rho$ + 0.8027$\omega$ &  2.58 &   3.0006$\rho$ + 0.8020$\omega$ \\ \hline

    2.6 &
    2.9986$\rho$ + 0.8014$\omega$ &  2.62 &   2.9967$\rho$ + 0.8008$\omega$ &  2.64 &
    2.9949$\rho$ + 0.8002$\omega$ \\ \hline

     2.66 &   2.9931$\rho$ + 0.7996$\omega$ &  2.68 &
    2.9915$\rho$ + 0.7991$\omega$ &  2.7 &   2.9899$\rho$ + 0.7986$\omega$ \\ \hline

      2.72 &
    2.9884$\rho$ + 0.7981$\omega$ &  2.74 &   2.9871$\rho$ + 0.7976$\omega$ &  2.76 &
    2.9859$\rho$ + 0.7972$\omega$ \\ \hline

    2.78 &   2.9848$\rho$ + 0.7968$\omega$ &  2.8 &
    2.9839$\rho$ + 0.7965$\omega$ &  2.82 &   2.9833$\rho$ + 0.7963$\omega$ \\ \hline

     2.84 &
    2.9828$\rho$ + 0.7961$\omega$ &  2.86 &   2.9826$\rho$ + 0.7959$\omega$ &  2.88 &
    2.9827$\rho$ + 0.7959$\omega$ \\ \hline

      2.9 &   2.9831$\rho$ + 0.7959$\omega$ &  2.92 &
    2.9840$\rho$ + 0.7961$\omega$ &  2.94 &   2.9854$\rho$ + 0.7964$\omega$ \\ \hline

     2.96 & 2.9876$\rho$ + 0.7970$\omega$ &  2.98 &   2.9909$\rho$ +
0.7978$\omega$ & 3. &   2.9973$\rho$ + 0.7995$\omega$ \\ \hline

\end{tabular}
\end{center}
\bigskip
\begin{center}
Table 4 \\
Values of the detection function $\lambda_{4}(h)$, for
$a=1/3, b=1/2, n=12.$ \\
\begin{tabular}{|c|c|c|c|c|c|}
  \hline
    $h$ & $\lambda_4(h)$ & $h$ & $\lambda_4(h)$ & $h$ & $\lambda_4(h)$ \\ \hline

    3&2.9973$\rho$+0.7995$\omega$ & 3.04 &
    3.0074$\rho$+0.8022$\omega$  &  3.08 &
    3.0118$\rho$+0.8032$\omega$ \\ \hline

     3.12 &
    3.0140$\rho$+0.8037$\omega$  &  3.16 &
    3.0149$\rho$+0.8037$\omega$ &  3.2 &
    3.0147$\rho$+0.8035$\omega$ \\ \hline

    3.24 &
    3.0138$\rho$+0.8030$\omega$ &   3.28 &
    3.0121$\rho$+0.8024$\omega$ &   3.32 &
    3.0098$\rho$+0.8016$\omega$ \\ \hline

    3.36 &
    3.0070$\rho$+0.8006$\omega$ &    3.4 &
    3.0038$\rho$+0.7995$\omega$ &    3.44 &
    3.0001$\rho$+0.7982$\omega$ \\ \hline

     3.48 &
    2.9960$\rho$+0.7969$\omega$ &   3.52 &
    2.9915$\rho$+0.7954$\omega$ &    3.56 &
    2.9867$\rho$+0.7939$\omega$ \\ \hline

     3.6 &
    2.9816$\rho$+0.7922$\omega$ &  3.64 &
    2.9762$\rho$+0.7905$\omega$ &  3.68 &
    2.9706$\rho$+0.7887$\omega$ \\ \hline

    3.72 &
    2.9646$\rho$+0.7868$\omega$ &  3.76 &
    2.9584$\rho$+0.7849$\omega$ &  3.8 &
    2.9520$\rho$+0.7829$\omega$ \\ \hline

     3.84 &
    2.9453$\rho$+0.7808$\omega$ &  3.88 &
    2.9384$\rho$+0.7787$\omega$ &  3.92 &
    2.9313$\rho$+0.7765$\omega$ \\ \hline

    3.96 &
    2.9240$\rho$+0.7742$\omega$ &  4. &  2.9165$\rho$+0.7719$\omega$
 &  4.2 &
    2.8764$\rho$+0.7596$\omega$ \\ \hline

     4.4 &  2.8322$\rho$+0.7461$\omega$ &  4.6 &
    2.7847$\rho$+0.7317$\omega$ &  4.8 &  2.7341$\rho$+0.7165$\omega$ \\ \hline

    5. &
    2.6808$\rho$+0.7005$\omega$ &  5.2 &  2.6250$\rho$+0.6838$\omega$ &  5.4 &
    2.5671$\rho$+0.6665$\omega$ \\ \hline

 5.6 &  2.5072$\rho$+0.6487$\omega$ &  5.8 &
    2.4455$\rho$+0.6304$\omega$ &  6. &  2.3820$\rho$+0.6116$\omega$ \\ \hline

     6.2 &
    2.3171$\rho$+0.5924$\omega$ &  6.4 &  2.2507$\rho$+0.5729$\omega$ &  6.6 &
    2.1830$\rho$+0.5530$\omega$ \\ \hline

    6.8 &  2.1141$\rho$+0.5327$\omega$ &  7. &
    2.0440$\rho$+0.5122$\omega$ &  7.2 &  1.9729$\rho$+0.4914$\omega$  \\ \hline

     7.4 &
    1.9008$\rho$+0.4704$\omega$ &  7.6 &  1.8278$\rho$+0.4491$\omega$ &  7.8 &
    1.7540$\rho$+0.4277$\omega$ \\ \hline

     8. &  1.6794$\rho$+0.4060$\omega$ &  8.2 &
    1.6041$\rho$+0.3842$\omega$ \\ \hline

\end{tabular}
\end{center}
\medskip

\medskip

\begin{figure}[h]\bce
\includegraphics{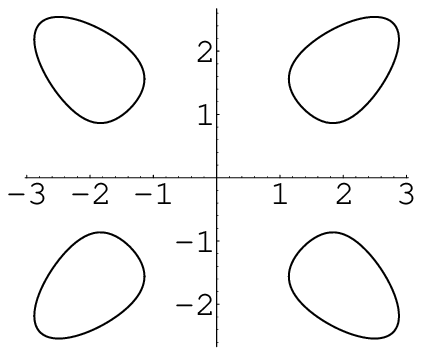}
\includegraphics{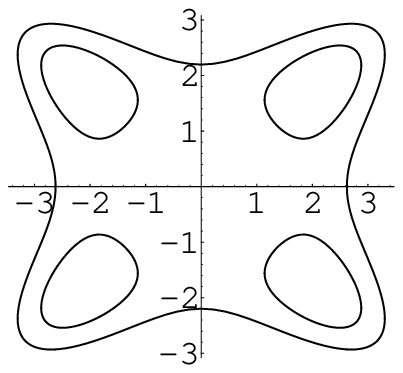}
\includegraphics{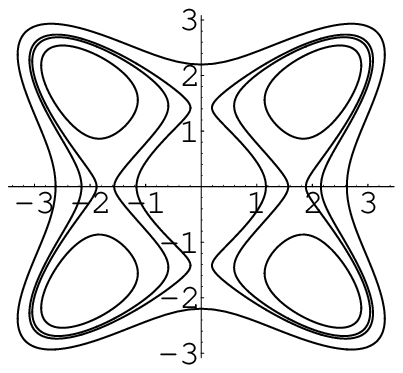}
\caption{ Distribution diagram corresponding to: a) four (left),
b) five (middle), c) nine (right) limit cycles of the system
(\ref{sis3n12})} \label{dist1} \ece
\end{figure}

\begin{figure}[h]\bce
\includegraphics{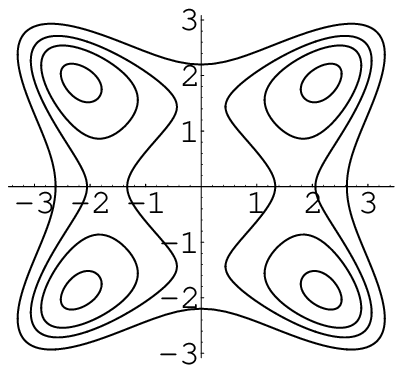}
\includegraphics{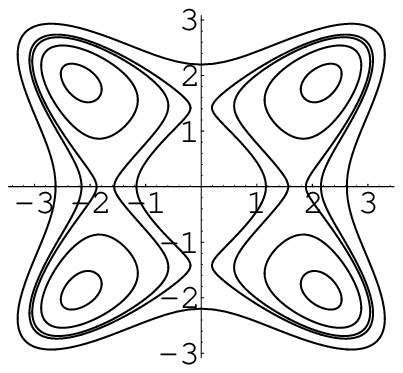}
\includegraphics{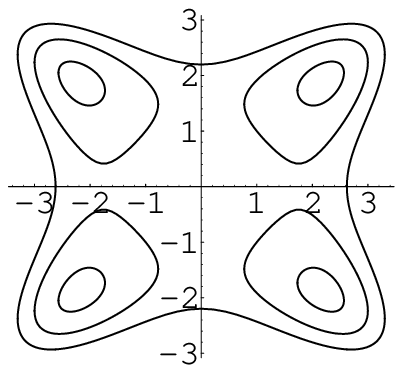}
\caption{ Distribution diagram corresponding to: a) eleven (left),
b) thirteen (middle), c) nine (right) limit cycles of the system
(\ref{sis3n12})} \label{dist2} \ece
\end{figure}

From the tables (1-4) we have the four detection functions, three
of which are shown in Fig.\ref{detfig1}, \ref{detfig2}. The values
of $\lambda_{2}$ are not plotted because they are too small in
comparison with the other values. From Proposition \ref{prop1},
and Fig.\ref{detfig1}, \ref{detfig2}, one gets:

\begin{theorem}
For $a=\frac{1}{3}, b=\frac{1}{2}, n=12, u =0.007$ and $v=-0.028$
and
$0<\varepsilon\ll 1$, we have the following distribution of the limit cycles: \\


 \noindent a) If $176.22<\lambda <242.6$,
the system (\ref{sis3n12}) has at least four limit cycles, one  in
the neighborhood of each orbit of type $\Gamma _{4}^{h}$, Fig.\ref{dist1}a),  \\
\\
b) If $242.6<\lambda <286.76$, the system (\ref{sis3n12}) has at
least five limit cycles, one  in the neighborhood of each orbit
of type $\Gamma _{4}^{h}$ and $\Gamma _{1}^{h}$, Fig.\ref{dist1}b), \\
\\
c) If $286.76<\lambda <288.49$, the system (\ref{sis3n12}) has at
least nine limit cycles, two  in the neighborhood of each orbit of
type $\Gamma _{3}^{h}$, and one in the neighborhood of
each orbit of type $\Gamma _{4}^{h}$, $\Gamma _{1}^{h}$, Fig.\ref{dist1}c), \\
\\
d) If $288.49<\lambda <288.92$,  the system (\ref{sis3n12}) has at
least eleven limit cycles, two  in the neighborhood of each orbit
of type $\Gamma _{4}^{h}$ and one  in the neighborhood of
each orbit of type $\Gamma _{3}^{h}$, $\Gamma _{1}^{h}$, Fig.\ref{dist2}a), \\
\\
e) If $288.92<\lambda < 289.99$, the system (\ref{sis3n12}) has at
least thirteen limit cycles, two  in the neighborhood of each
orbit of type $\Gamma _{4}^{h}$, $\Gamma _{3}^{h}$ and one
in the neighborhood of the orbit of type $\Gamma _{1}^{h}$, Fig.\ref{dist2}b), \\
\\
f) If $289.99<\lambda <290.82$, the system (\ref{sis3n12}) has at
least nine limit cycles, two in the neighborhood of each orbit of
type $\Gamma _{4}^{h}$ and one in the neighborhood of the orbit of
type $\Gamma _{1}^{h}$, Fig.\ref{dist2}c) .
\end{theorem}

 From Fig.\ref{detfig1}, \ref{detfig2}
we could obtain more results but we listed only the important
cases.

\section{Conclusion}

The system (\ref{sis3n12}) corresponding to $a=\frac{1}{3},
b=\frac{1}{2}, n=12, u =0.007, v=-0.028, 0<\varepsilon\ll 1$, and
$288.49<\lambda <288.92$, has at least eleven limit cycles while
for $288.92<\lambda < 289.99$ it has at least thirteen limit
cycles. The Abelian integral method was used. Through numerical
explorations we have drawn the shape of the graphs of the
detection functions, from which we determined the number and
distribution of limit cycles. A natural question is: Is it
possible to obtain more limit cycles for higher perturbations?

\section{Acknowledgements}

This work was partially supported through a European Community
Marie Curie Fellowship, and in the framework of the CTS, contract
number HPMT-CT-2001-00278.

\end{document}